# Local asymptotic normality for Cox–Ingersoll–Ross drift parameters under high-frequency sampling and the strict Feller condition

Boyuan Ning*

*Department of Pure and Applied Mathematics, Waseda University, 3-4-1 Ohkubo, Shinjuku-ku, Tokyo 169-8555, Japan*

**Abstract**

For the ergodic Cox–Ingersoll–Ross (CIR) diffusion with known diffusion scale and a fixed initial value, we prove local asymptotic normality (LAN) of the exactly sampled experiment for the drift parameters whenever the Feller index exceeds 1, assuming only a vanishing mesh and a diverging observation horizon. This improves the previously known condition that the index exceed $(11+\sqrt{89})/2$. The proof cuts off the continuous score near 0, recovers the remainder by exact conditional centering on each sampling interval, and transfers LAN by a conditional-likelihood argument.

*Keywords:* Cox–Ingersoll–Ross process, local asymptotic normality, high-frequency observations, drift estimation, score recovery, information loss

*2020 MSC:* Primary 62M05; Secondary 62F12, 60J60.

*Corresponding author

*Email address:* ningboyuan@aoni.waseda.jp (Boyuan Ning)

## 1. Experiment and main statement

For $\theta = (a, b) \in \Theta = (\sigma, \infty) \times (0, \infty)$, let $\mathbb{P}_\theta$ denote the law of the solution to

$$dX_t = (a - bX_t)dt + \sqrt{2\sigma X_t}dW_t, \qquad X_0 = x_0, \tag{1.1}$$

with fixed known $\sigma > 0$ and $x_0 > 0$. We work on the canonical continuous-path space with the raw natural filtration $\mathcal{F}_t = \sigma(X_s : 0 \le s \le t)$. Likelihood ratios are formed on these finite-time $\sigma$-fields. Write $r = a/\sigma$ for the Feller index. We observe $(X_{k\Delta_n})_{k=0}^n$, with

$$\Delta_n \longrightarrow 0, \qquad T_n = n\Delta_n \longrightarrow \infty. \tag{1.2}$$

Let $\mathcal{G}_n = \sigma(X_{k\Delta_n} : 0 \le k \le n)$, and denote by $Q_{\theta,n}$ the distribution of this observation vector. The family $\{Q_{\theta,n} : \theta \in \Theta\}$ defines the exact sampled experiment.

Fix $\theta_0 = (a_0, b_0) \in \Theta$, and abbreviate $\mathbb{P}_0 = \mathbb{P}_{\theta_0}$ and $\mathbb{E}_0 = \mathbb{E}_{\theta_0}$. Under this law, define

$$M_T = (M_T^a, M_T^b)^\top := \Big( \int_0^T \frac{dW_s}{\sqrt{2\sigma X_s}}, -\int_0^T \sqrt{\frac{X_s}{2\sigma}} dW_s \Big)^\top, \tag{1.3}$$

and $\widehat{M}_n := \mathbb{E}_0[M_{T_n} \mid \mathcal{G}_n]$. By Itô's formula and the defining equation of $X$,

$$\begin{aligned} M_T^a &= \frac{1}{2\sigma}\Big[\log \frac{X_T}{X_0} - (a_0 - \sigma)\int_0^T \frac{ds}{X_s} + b_0 T\Big], \\ M_T^b &= -\frac{1}{2\sigma}\Big[X_T - X_0 - a_0 T + b_0 \int_0^T X_s ds\Big]. \end{aligned} \tag{1.4}$$

Thus both components are measurable functions of the observed continuous path. By Lemma 2.1, $M_T$ is the continuous-observation score at $\theta_0$.

**Theorem 1.1** (Main result). *Suppose $a_0 > \sigma$, $b_0 > 0$, and* (1.2) *holds. For every deterministic sequence $h_n \to h \in \mathbb{R}^2$, the parameters $\theta_n = \theta_0 + h_n/\sqrt{T_n}$ belong to $\Theta$ for all sufficiently large $n$, and*

$$\log \frac{dQ_{\theta_n,n}}{dQ_{\theta_0,n}}(X_0,\dots,X_{T_n}) = h_n^\top \frac{\widehat{M}_n}{\sqrt{T_n}} - \frac{1}{2} h_n^\top I_0 h_n + o_{\mathbb{P}_0}(1), \qquad \frac{\widehat{M}_n}{\sqrt{T_n}} \xrightarrow{\mathcal{L}} \mathcal{N}_2(0, I_0), \tag{1.5}$$

*where*

$$I_0 = \frac{1}{2\sigma}\begin{pmatrix} b_0/(a_0-\sigma) & -1 \\ -1 & a_0/b_0 \end{pmatrix} \tag{1.6}$$

*is positive definite. In particular, the discrete and continuous experiments have the same LAN rate and limiting information. No additional relation between $T_n$ and $\Delta_n$ is assumed. At $r = 1$, $\pi_{\theta_0}(x^{-1}) = \infty$, so the continuous-observation Fisher information in the $a$-direction is infinite.*

*Relation to earlier work.* Ben Alaya et al. (2020, Theorem 2.1) established LAN for this experiment under the same sampling conditions (1.2) and a fixed initial value, assuming $r > (11+\sqrt{89})/2 \approx 10.22$. Theorem 1.1 extends this result to the full finite-information range $r > 1$. The transfer of LAN from the continuous to the sampled experiment follows classical arguments (Le Cam and Yang, 1988), and recovering a continuous-observation central sequence from discrete data is also used by Becheri et al. (2016). The model-specific step is the score recovery bound of Proposition 3.2, which requires an inverse moment of order only slightly above one. A cutoff near zero and exact interval centering yield orthogonal recovery errors (Section 3), allowing

LAN to transfer through conditional likelihoods (Section 4).

Under stationary initialization, the occupation-time bounds of Altmeyer and Chorowski (2018), together with (1.4), yield an $O(T\Delta)$ mean-square recovery bound for the $a$-score when $r > 2$. The cutoff in Proposition 3.2 also handles $1 < r \leq 2$, where $x^{-1} \notin L^2(\pi_{\theta_0})$. Detailed comparisons are given in Appendix B of the online appendix.

The extended range also matters for estimation. The asymptotic covariance of the discretized estimator of Chernova et al. (2024, Theorem 2.2), obtained for $r > 2$ under a polynomial sampling scheme satisfying $T_n\Delta_n = n\Delta_n^2 \to 0$, and the drift block of that of the Gaussian quasi-likelihood estimator of Cheng et al. (2024), obtained under stationary initialization and $r > 5$, both equal $I_0^{-1}$. Cheng et al. (2024, Remark 2.4) obtain efficiency from the LAN result of Ben Alaya et al. (2020), hence only for $r > (11 + \sqrt{89})/2$, and suggest that it should hold more generally; Theorem 1.1 identifies $I_0^{-1}$ as the asymptotic covariance lower bound for regular drift estimators throughout $r > 1$.

## 2. CIR inputs and the continuous experiment

For $\theta \in \Theta$, the solution of (1.1) is non-explosive and strictly positive, and its invariant distribution is $\pi_\theta(dx) = \frac{(b/\sigma)^r}{\Gamma(r)} x^{r-1} e^{-bx/\sigma} dx$, $x > 0$. We use the standard fixed-initial-value moment bounds

$$\sup_{t\geq 0} \mathbb{E}_\theta[X_t^q] < \infty \quad (q > 0), \qquad \sup_{t\geq 0} \mathbb{E}_\theta[X_t^{-p}] < \infty \quad (0 < p < r). \tag{2.1}$$

For a fixed starting point, both bounds follow from Ben Alaya and Kebaier (2013, Proposition 3(2)); the inverse-moment bound is also recorded in Ben Alaya et al. (2020, Lemma 3.1(i)). We also use the CIR ergodic theorem, from a fixed positive initial value, for $\pi_\theta$-integrable functions; see Ben Alaya and Kebaier (2013, proof of Theorem 5) and Ben Alaya et al. (2020, Equation (2.2)). In particular,

$$\frac{1}{T}\int_0^T X_s ds \longrightarrow \frac{a}{b}, \qquad \frac{1}{T}\int_0^T \frac{ds}{X_s} \longrightarrow \frac{b}{a-\sigma} \quad \mathbb{P}_\theta\text{-a.s.} \tag{2.2}$$

The continuous-observation LAN argument below is classical; see Overbeck (1998). The proof of the following finite-horizon likelihood lemma is given in Appendix A.1 of the online appendix.

**Lemma 2.1** (Finite-horizon likelihood)**.** *For any $\theta, \bar\theta \in \Theta$ and $T < \infty$, the path laws restricted to $\mathcal{F}_T$ are equivalent. Under $\mathbb{P}_\theta$, their likelihood ratio is*

$$L_T(\bar\theta, \theta) = \exp\Big\{(\bar\theta - \theta)^\top M_T^\theta - \frac{1}{2}(\bar\theta - \theta)^\top J_T(\bar\theta - \theta)\Big\},$$

*where $M_T^\theta$ is defined as in* (1.3) *under $\mathbb{P}_\theta$ and $J_T = \frac{1}{2\sigma}\int_0^T \begin{pmatrix} X_s^{-1} & -1 \\ -1 & X_s \end{pmatrix} ds$. In particular, $\mathbb{E}_\theta[L_T(\bar\theta, \theta)] = 1$.*

For deterministic $h_n \to h$ and $\theta_n = \theta_0 + h_n/\sqrt{T_n}$, Lemma 2.1 gives the *exact* expansion

$$\log L_{T_n}(\theta_n, \theta_0) = h_n^\top Z_n - \frac{1}{2} h_n^\top \frac{J_{T_n}}{T_n} h_n, \qquad Z_n = \frac{M_{T_n}}{\sqrt{T_n}}.$$

The moment bounds give $\mathbb{E}_0[\operatorname{tr} J_t] < \infty$ for every finite $t$, so $M$ is a continuous square-integrable martingale with bracket $J$. For $\mathcal{M}_n(u) = M_{T_n u}/\sqrt{T_n}$ and every $u \geq 0$,

$$\langle \mathcal{M}_n \rangle_u = J_{T_n u}/T_n \longrightarrow u I_0 \quad \mathbb{P}_0\text{-a.s.},$$

by (2.2). All jumps of $\mathcal{M}_n$ and its bracket vanish. The martingale CLT (e.g. Whitt, 2007, Theorem 2.1(ii)) therefore gives $Z_n \xrightarrow{\mathcal{L}} \mathcal{N}_2(0, I_0)$. Since $(h_n)$ is bounded, $h_n^\top (J_{T_n}/T_n - I_0) h_n \xrightarrow{\mathbb{P}_0} 0$, proving continuous LAN along these sequences. The condition $a_0 > \sigma$ enters through $\pi_{\theta_0}(x^{-1}) < \infty$ in the $a$-direction.

## 3. Recovery of the continuous score

Throughout this section the parameter is $\theta_0$; write $a = a_0$, $b = b_0$. For $T = n\Delta$, set $t_k = k\Delta$ and $\mathcal{G}_{n,\Delta} = \sigma(X_{t_k} : 0 \leq k \leq n)$. Constants in the recovery bounds may depend on $(a, b, \sigma, x_0)$, on a fixed cutoff function, and on $\eta$ when indicated, but not on $T, n, \Delta, \varepsilon$. Dependence on $\varepsilon$ is displayed explicitly. Write $K_p^- := \sup_{t \geq 0} \mathbb{E}_0[X_t^{-p}] < \infty \quad (0 < p < r)$ and $K_2^+ := \sup_{t \geq 0} \mathbb{E}_0[X_t^2] < \infty$. Let $\mathsf{S}_t g(x) = \mathbb{E}_{\theta_0, x}[g(X_t)]$ denote the Markov semigroup, where $x$ denotes the initial state, and put the infinitesimal generator $\mathcal{L}g(x) = (a - bx)g'(x) + \sigma x g''(x)$.

**Lemma 3.1** (Exact interval centering)**.** *Let $G \in C^2((0,\infty))$, put $f = \mathcal{L}G$, and assume, for each finite $T$, that $G(X_t) \in L^2(\mathbb{P}_0)$ for $t \leq T$ and*

$$\int_0^T \mathbb{E}_0\big[2\sigma X_s |G'(X_s)|^2 + |f(X_s)|^2\big] ds < \infty.$$

*Suppose the Itô decomposition*

$$G(X_t) - G(X_0) = M_t^G + \int_0^t f(X_s)ds, \qquad M_t^G = \int_0^t \sqrt{2\sigma X_s} G'(X_s) dW_s$$

*holds. Its martingale term is then a true square-integrable martingale. Define the $\mathcal{G}_{n,\Delta}$-measurable statistic*

$$A_{n,G} := \sum_{k=0}^{n-1} \{G(X_{t_{k+1}}) - \mathsf{S}_\Delta G(X_{t_k})\}. \tag{3.1}$$

*Then*

$$\mathbb{E}_0\big[|A_{n,G} - M_T^G|^2\big] \le \Delta \int_0^T \mathbb{E}_0\big[|f(X_s)|^2\big] ds.$$

*Proof.* Let $B_k = \int_{t_k}^{t_{k+1}} f(X_s)ds$. The martingale increment has conditional mean zero. The Markov property therefore gives

$$\mathsf{S}_\Delta G(X_{t_k}) - G(X_{t_k}) = \mathbb{E}_0[B_k \mid \mathcal{F}_{t_k}],$$

and consequently

$$A_{n,G} - M_T^G = \sum_{k=0}^{n-1} \{B_k - \mathbb{E}_0[B_k \mid \mathcal{F}_{t_k}]\}.$$

Write each summand as $R_k$. It is square integrable, $\mathcal{F}_{t_{k+1}}$-measurable, and has conditional mean zero given $\mathcal{F}_{t_k}$. For $j < k$, $R_j$ is $\mathcal{F}_{t_k}$-measurable and $\mathbb{E}_0[R_j R_k] = \mathbb{E}_0[R_j \mathbb{E}_0(R_k \mid \mathcal{F}_{t_k})] = 0$. Thus the errors are orthogonal.

Conditional variance contraction and Cauchy–Schwarz give

$$\mathbb{E}_0\big[|A_{n,G} - M_T^G|^2\big] \le \sum_{k=0}^{n-1} \mathbb{E}_0[B_k^2] \le \Delta \sum_{k=0}^{n-1} \int_{t_k}^{t_{k+1}} \mathbb{E}_0\big[|f(X_s)|^2\big] ds. \qquad \square$$

**Proposition 3.2** (Score recovery bound)**.** *For every fixed* $\eta \in (0, r-1)$ *there is* $C_\eta < \infty$ *such that, whenever* $0 < \Delta \le 1$ *and* $T = n\Delta$,

$$\frac{1}{T}\mathbb{E}_0\big[\|M_T - \mathbb{E}_0[M_T \mid \mathcal{G}_{n,\Delta}]\|^2\big] \le C_\eta(\Delta^{\eta/(\eta+2)} + \Delta). \tag{3.2}$$

*In particular, under* (1.2), $\frac{1}{\sqrt{T_n}}(M_{T_n} - \widehat{M}_n) \to 0$ *in* $L^2(\mathbb{P}_0)$.

*Proof. The singular a-direction.* Choose a smooth function $\chi : [0,\infty) \to [0,1]$ which vanishes on $[0,1/2]$ and equals one on $[1,\infty)$. For $0 < \varepsilon \le 1$, set $\chi_\varepsilon(x) = \chi(x/\varepsilon)$ and

$$M_T^{a,\varepsilon} := \int_0^T \frac{\chi_\varepsilon(X_s)}{\sqrt{2\sigma X_s}} dW_s.$$

Since $1 + \eta < r$, the moment bound (2.1) and Itô isometry imply

$$\begin{aligned}\mathbb{E}_0\big[|M_T^a - M_T^{a,\varepsilon}|^2\big] &\le \frac{1}{2\sigma}\int_0^T \mathbb{E}_0[X_s^{-1}\mathbf{1}_{\{X_s \le \varepsilon\}}]ds \le \frac{\varepsilon^\eta}{2\sigma}\int_0^T \mathbb{E}_0[X_s^{-(1+\eta)}]ds \\ &\le \frac{K_{1+\eta}^-}{2\sigma} T\varepsilon^\eta. \end{aligned} \tag{3.3}$$

Now take the primitive

$$G_\varepsilon(x) = \int_1^x \frac{\chi_\varepsilon(y)}{2\sigma y} dy.$$

Its martingale part is exactly $M^{a,\varepsilon}$, since $\sqrt{2\sigma x}G_\varepsilon'(x) = \chi_\varepsilon(x)/\sqrt{2\sigma x}$. Direct

differentiation gives

$$f_\varepsilon(x) = \mathcal{L}G_\varepsilon(x) = \chi_\varepsilon(x)\Big(\frac{a-\sigma}{2\sigma x} - \frac{b}{2\sigma}\Big) + \frac{1}{2}\chi_\varepsilon'(x), \qquad \|f_\varepsilon\|_\infty \le C_f\varepsilon^{-1},$$

where $C_f < \infty$ is independent of $\varepsilon$. Also, with $H_\varepsilon(x) = \sqrt{2\sigma x}G_\varepsilon'(x)$,

$$\|H_\varepsilon\|_\infty \le (\sigma\varepsilon)^{-1/2}, \qquad |G_\varepsilon(x)| \le C(1 + |\log\varepsilon| + \log^+ x).$$

Hence the endpoint values of $G_\varepsilon(X_t)$ are square integrable. For fixed $\varepsilon, T$, apply Itô's formula up to $\tau_m = \inf\{t : X_t \notin (1/m, m)\}$. Positivity and non-explosion give $\tau_m \uparrow \infty$ almost surely. Boundedness of $H_\varepsilon$ and $f_\varepsilon$ makes both stopped integrals converge in $L^2$, by Itô isometry and dominated convergence. The stopped endpoint converges almost surely to $G_\varepsilon(X_t)$; the stopped identity therefore identifies its $L^2$ limit and proves the unstopped decomposition. This verifies all hypotheses of Lemma 3.1, which yields

$$\mathbb{E}_0\big[|A_{n,G_\varepsilon} - M_T^{a,\varepsilon}|^2\big] \le C_f^2 T\Delta\varepsilon^{-2}. \tag{3.4}$$

Since $A_{n,G_\varepsilon} \in L^2(\mathcal{G}_{n,\Delta})$, the projection property of conditional expectation, together with (3.3) and (3.4), gives

$$\mathbb{E}_0\big[|M_T^a - \mathbb{E}_0[M_T^a \mid \mathcal{G}_{n,\Delta}]|^2\big] \le \mathbb{E}_0\big[|M_T^a - A_{n,G_\varepsilon}|^2\big] \le T\Big(\frac{K_{1+\eta}^-}{\sigma}\varepsilon^\eta + 2C_f^2\Delta\varepsilon^{-2}\Big). \tag{3.5}$$

This bound holds for every $0 < \varepsilon, \Delta \le 1$ and $T = n\Delta$, with no other dependence on $\varepsilon$ or $T$. Taking $\varepsilon = \Delta^{1/(\eta+2)}$ therefore bounds it by $(K_{1+\eta}^-/\sigma + 2C_f^2)T\Delta^{\eta/(\eta+2)}$.

*The b-direction.* Take $G_b(x) = -x/(2\sigma)$. Its martingale part is $M^b$, and $\mathcal{L}G_b(x) = -a/(2\sigma) + bx/(2\sigma)$. By the uniform second moment in (2.1), Lemma 3.1 gives

$$\mathbb{E}_0\big[|M_T^b - \mathbb{E}_0[M_T^b \mid \mathcal{G}_{n,\Delta}]|^2\big] \le \mathbb{E}_0\big[|M_T^b - A_{n,G_b}|^2\big] \le \frac{a^2 + b^2 K_2^+}{2\sigma^2} T\Delta, \quad (3.6)$$

The same moment bounds make the endpoints, stochastic integral, and drift integral square integrable. Adding (3.5) with the selected cutoff and (3.6) proves (3.2). □

## 4. Transfer from path likelihoods to sampled likelihoods

The following LAN-preservation lemma follows classical arguments (Le Cam and Yang, 1988); see also Davies (1985, Theorem 4.1). A self-contained proof is given in Appendix A.2 of the online appendix.

**Lemma 4.1** (Conditional-likelihood transfer)**.** *Let $h_n \to h \in \mathbb{R}^d$ be deterministic. For each $n$, let $P_{n,h_n} \ll P_{n,0}$ be probability measures on a full experiment, with likelihood ratio $L_n = L_n(h_n) > 0$. Suppose that*

$$\log L_n = h_n^\top Z_n - \frac{1}{2} h_n^\top I h_n + o_{P_{n,0}}(1), \qquad Z_n \xrightarrow{\mathcal{L}} \mathcal{N}_d(0, I),$$

*where $I$ is deterministic, symmetric, and non-negative definite. Let $\mathcal{G}_n$ be a sub-$\sigma$-field and let $U_n$ be $\mathcal{G}_n$-measurable with $U_n - Z_n \to 0$ in $P_{n,0}$-probability. Then the restricted likelihood ratio $\widetilde{L}_n = E_{n,0}[L_n \mid \mathcal{G}_n]$ satisfies*

$$\log \widetilde{L}_n = h_n^\top U_n - \frac{1}{2} h_n^\top I h_n + o_{P_{n,0}}(1).$$

*Proof of Theorem 1.1.* Take the full experiment to be observation of the path on $[0, T_n]$ and the reduced experiment to be its sampled coordinates. By Lemma 2.1, the full likelihood ratio is a probability density with expectation one. The Radon–Nikodým derivative on the reduced experiment is exactly

$$\frac{dQ_{\theta_n,n}}{dQ_{\theta_0,n}}(X_0, \ldots, X_{T_n}) = \mathbb{E}_0[L_{T_n}(\theta_n, \theta_0) \mid \mathcal{G}_n]. \tag{4.1}$$

Indeed, both sides of (4.1), when integrated over a grid-measurable event under $\mathbb{P}_0$, give its probability under $\mathbb{P}_{\theta_n}$. Section 2 proves continuous LAN along every deterministic $h_n \to h$, with $Z_n = M_{T_n}/\sqrt{T_n}$ and information $I_0$. Proposition 3.2 shows that the observable vector $U_n = \widehat{M}_n/\sqrt{T_n}$ satisfies $U_n - Z_n \to 0$ in $L^2$. Lemma 4.1 now gives the expansion (1.5), while Slutsky's theorem gives $U_n \xrightarrow{\mathcal{L}} \mathcal{N}_2(0, I_0)$. Since $\theta_0$ is an interior point and $(h_n)$ is bounded, all local parameters lie in $\Theta$ for sufficiently large $n$. □

**Funding**

This research did not receive any specific grant from funding agencies in the public, commercial, or not-for-profit sectors.

**Data availability**

No datasets were generated or analyzed in this theoretical study.

## Declaration of generative AI and AI-assisted technologies in the manuscript preparation process

During the preparation of this work, the author used ChatGPT and Codex (OpenAI) and Claude (Anthropic) to assist with literature searches, exploration and assessment of mathematical arguments, and manuscript drafting and revision. After using these tools, the author reviewed and edited the content as needed and takes full responsibility for the content of the published article.

## Acknowledgements

The author thanks Professor Yasutaka Shimizu at Waseda University for his guidance, helpful discussions, and valuable suggestions.

# Online Appendix to

Local asymptotic normality for Cox–Ingersoll–Ross drift parameters under high-frequency sampling and the strict Feller condition

Boyuan Ning

Department of Pure and Applied Mathematics, Waseda University
Tokyo 169-8555, Japan

ningboyuan@aoni.waseda.jp

This appendix provides the proofs of Lemmas 2.1 and 4.1 of the main article and detailed comparisons with related work. Notation follows the main article. Unless an external source is specified, references to numbered results and equations without an appendix prefix refer to the main article.

# A Technical Proofs

## A.1 Finite-horizon likelihood

*Proof of Lemma 2.1 of the main article.* Set $\lambda(x) = [(\bar{a} - a) - (\bar{b} - b)x]/\sqrt{2\sigma x}$ and let $\tau_m = \inf\{t : X_t \notin (1/m, m)\}$, with $m$ large enough that $x_0 \in (1/m, m)$. Write $L_T^{(m)}$ for the stochastic exponential with integrand $\mathbf{1}_{\{s \le \tau_m\}}\lambda(X_s)$. This integrand is bounded, so $\mathbb{E}_\theta[L_T^{(m)}] = 1$ for each finite $T$. Define $\widetilde{\mathbb{P}}_m$ on $\mathcal{F}_T$ by $d\widetilde{\mathbb{P}}_m = L_T^{(m)} d\mathbb{P}_\theta$. Girsanov's theorem gives drift $a - bx + \sqrt{2\sigma x}\lambda(x) = \bar{a} - \bar{b}x$ before $\tau_m$. The coefficients are Lipschitz on the stopped interval, so uniqueness identifies the law of the stopped path under $\widetilde{\mathbb{P}}_m$ with that under $\mathbb{P}_{\bar{\theta}}$.

For $A \in \mathcal{F}_T$, the event $A \cap \{\tau_m > T\}$ is measurable with respect to $\mathcal{F}_{T \wedge \tau_m}$, and $L_T^{(m)} = L_T$ on $\{\tau_m > T\}$. Hence

$$\mathbb{E}_\theta[L_T(\bar{\theta}, \theta)\mathbf{1}_{A \cap \{\tau_m > T\}}] = \widetilde{\mathbb{P}}_m(A \cap \{\tau_m > T\}) = \mathbb{P}_{\bar{\theta}}(A \cap \{\tau_m > T\}).$$

Under $\mathbb{P}_\theta$ and $\mathbb{P}_{\bar{\theta}}$, positivity and non-explosion imply $\tau_m \uparrow \infty$ almost surely. Monotone convergence on both sides proves $\mathbb{E}_\theta[L_T \mathbf{1}_A] = \mathbb{P}_{\bar{\theta}}(A)$, including $\mathbb{E}_\theta[L_T] = 1$. On every finite horizon, a positive continuous path stays in a compact subset of $(0, \infty)$, so the unstopped exponent is finite and $L_T > 0$. This proves equivalence and true martingality on finite horizons. □

## A.2 Conditional-likelihood transfer

*Proof of Lemma 4.1 of the main article.* Write $R_n = h_n^\top U_n - h_n^\top I h_n/2$. Boundedness of $(h_n)$ gives $\log L_n - R_n \to 0$ in probability and

$$L_n \xrightarrow{\mathcal{L}} L = \exp\{h^\top Z - h^\top I h/2\}, \qquad Z \sim \mathcal{N}_d(0, I), \qquad \mathbb{E}[L] = 1.$$

Because $E_{n,0}[L_n] = 1$ and the limit has mean one, the likelihood ratios are uniformly integrable. More explicitly, for every $K > 0$,

$$E_{n,0}[(L_n - K)_+] = 1 - E_{n,0}[(L_n \wedge K)] \longrightarrow 1 - \mathbb{E}[(L \wedge K)] = \mathbb{E}[(L - K)_+]. \tag{A.1}$$

The right side tends to zero as $K \to \infty$. Using $E_{n,0}[L_n \mathbf{1}_{\{L_n > 2K\}}] \le 2E_{n,0}[(L_n - K)_+]$ and accounting for finitely many initial indices proves uniform integrability.

The function $x \mapsto e^x \wedge K$ is bounded and Lipschitz. Hence $\log L_n - R_n \to 0$ in probability implies

$$E_{n,0}\big[|(L_n \wedge K) - (e^{R_n} \wedge K)|\big] \longrightarrow 0.$$

Combining this with (A.1) gives

$$\limsup_{n\to\infty} E_{n,0}\big[|L_n - (e^{R_n}\wedge K)|\big] \le \mathbb{E}[(L-K)_+].$$

The capped variable is $\mathcal{G}_n$-measurable. Conditional expectation is an $L^1$ contraction, so the same upper bound applies with $L_n$ replaced by $\widetilde{L}_n$. For every $\delta > 0$,

$$P_{n,0}(|\widetilde{L}_n - e^{R_n}| > \delta) \le \delta^{-1} E_{n,0}\big[|\widetilde{L}_n - (e^{R_n}\wedge K)|\big] + P_{n,0}(e^{R_n} > K).$$

Since $e^{R_n} \xrightarrow{\mathcal{L}} L$, letting first $n\to\infty$ and then $K\to\infty$ (through continuity points) proves $\widetilde{L}_n - e^{R_n} \to 0$ in probability. Finally, $R_n$ is tight, so $e^{-R_n}$ is tight. It follows that $\widetilde{L}_n/e^{R_n}\to 1$ in probability. Since $L_n > 0$ almost surely and $E_{n,0}[L_n] = 1$, its conditional expectation is finite and strictly positive almost surely. Continuity of the logarithm at one proves the assertion. □

# B Related literature

**Comparison with the Malliavin proof for CIR.** The restriction $r > (11+\sqrt{89})/2$ in Ben Alaya et al. (2020) is a sufficient condition arising from a particular proof, not a necessary condition for LAN. Their Remark 2.5 explicitly attributes it to the Hölder estimates used to control the remainder terms in Lemmas 5.3 and 5.6. Proposition 4.4 first gives a score representation under $r > (9+\sqrt{57})/2$, and then a decomposition of the Skorohod weights under $r > 5+3\sqrt{2}$; the final remainder estimates impose the stronger threshold above. These conditions accompany estimates of stochastic-flow derivatives and Malliavin remainder terms.

The present proof instead recovers the continuous martingale score in $L^2$, without expanding a differentiated transition likelihood. Its singular tail is controlled by $\sup_{t\ge 0}\mathbb{E}_0[X_t^{-(1+\eta)}] < \infty$ for any fixed $0 < \eta < r-1$. For the retained part, the bound $\|\mathcal{L}G_\varepsilon\|_\infty \le C\varepsilon^{-1}$ and exact interval centering give orthogonal errors. Equation (3.5) therefore bounds the normalized loss by $C_\eta(\varepsilon^\eta + \Delta\varepsilon^{-2})$, uniformly in $T$. Only its convergence to zero is needed for LAN. This explains the weaker moment requirement and the short proof. The comparison concerns the ergodic case: Ben Alaya et al. (2020) also treat critical and supercritical regimes.

**Why the cutoff matters.** For the uncut $a$-score, take $G(x) = \log x/(2\sigma)$. Then

$$\mathcal{L}G(x) = \frac{a_0-\sigma}{2\sigma x} - \frac{b_0}{2\sigma}.$$

If $r = a_0/\sigma > 2$, the uniform inverse second-moment bound in (2.1) allows Lemma 3.1 to be applied directly. Together with the $b$-direction estimate, this gives

$$\mathbb{E}_0\big[\|M_T - \mathbb{E}_0[M_T \mid \mathcal{G}_{n,\Delta}]\|^2\big] \le CT\Delta$$

from the fixed initial value used here.

There is also a direct connection with existing occupation-time approximation results. Under stationary initialization, $x^{-1}\in L^2(\pi_{\theta_0})$ if and only if $r > 2$. The CIR process is reversible, so Altmeyer and Chorowski (2018, Theorem 2.2, $s=0$) give

$$\mathbb{E}_{\pi_{\theta_0}}\Big[\Big|\int_0^T X_s^{-1}ds - \Delta\sum_{k=0}^{n-1} X_{k\Delta}^{-1}\Big|^2\Big] \le CT\Delta.$$

The endpoint representation (1.4) and the projection property of conditional expectation yield the corresponding $a$-score recovery bound.

For $1 < r \le 2$, the direct application fails because $x^{-1}\notin L^2(\pi_{\theta_0})$. Proposition 3.2 handles this range by cutting off the martingale score, controlling its discarded part with an inverse moment of order $1+\eta < r$, and applying exact interval centering to the retained part. The resulting rate suffices for LAN.

**General diffusion results.** Gobet (2002) and Tran and Ngô (2022) establish discrete-observation local asymptotic results for regular diffusions. In particular, Section 5.1.1 of Tran and Ngô (2022) explicitly includes ergodic drift LAN. Their global Lipschitz and uniform ellipticity assumptions exclude the square-root coefficient at zero. The Lamperti transform $Y = \sqrt{2X/\sigma}$ gives unit diffusion but drift $(r - 1/2)/Y - bY/2$, whose derivative is unbounded near zero. Hue et al. (2026) study LAQ at critical parameter points under the same global Lipschitz and uniform ellipticity assumptions.

Sørensen (2024) develops a general theory of efficient high-frequency estimation based on approximate martingale estimating functions. Condition 2.2(4) requires globally Lipschitz drift and diffusion coefficients, excluding the square-root coefficient. Genon-Catalot and Larédo (2014) establish global asymptotic equivalence results relating continuous diffusion observations, discrete observations, and Euler scheme experiments. Their assumptions include a known diffusion coefficient bounded away from zero and the sampling condition $T_n\Delta_n = n\Delta_n^2 \to 0$.

The recovery-and-transfer strategy is also used by Becheri et al. (2016, Proposition 2). Their running Example 1 is a CIR process with Hawkes jumps under $r > 1$. Their discrete-observation result, however, assumes a stationary start, $T_n\Delta_n \to 0$ (Assumption 8), a jump intensity bounded away from zero, and a volatility bounded away from zero (Assumption 9(iii)); the last condition is used in their Appendix A.2 to control the step-process approximation of the drift-score weight, which for the CIR drift parameters is $\nabla_\theta(a - bx)/(2\sigma x) = (x^{-1}, -1)^\top/(2\sigma)$. This lower-bound condition fails for $\sqrt{2\sigma x}$ near zero, so their result does not cover (1.1).

For jump-type CIR processes, Ben Alaya et al. (2025) establish local asymptotic results under continuous and high-frequency discrete observations for the growth-rate parameter $b$, with the immigration parameter $a$ and diffusion coefficient known. Their ergodic LAN result concerns the $b$-direction; joint drift inference in the present paper additionally requires recovery of the singular $a$-score.

**Other CIR observation schemes and estimators.** Overbeck and Rydén (1997) study CIR inference at a fixed positive sampling interval. More recently, Amengual and Puig-Codina (2026, Theorem 2) obtain exact-likelihood LAN jointly for the drift parameters and the diffusion scale, from a fixed initial value at fixed $\Delta > 0$, without imposing the Feller condition. Theorem 1.1 treats the high-frequency regime $\Delta_n \to 0$ and $T_n \to \infty$, with limiting information equal to that of the continuous-observation experiment.

Cheng et al. (2024) study high-frequency Gaussian quasi-likelihood estimation under stationary initialization, jointly estimating the drift and diffusion parameters without an additional sampling coupling. Their parameter condition is $r > 5$. Their Remark 2.4 discusses the possibility of extending the efficiency conclusion beyond the parameter restriction in Ben Alaya et al. (2020).

Chernova et al. (2024, Theorem 2.2) study discretized versions of the continuous-time maximum likelihood estimator and assume $r > 2$; their Remark 2.3 attributes this restriction to the analysis of Riemann sums of $X^{-1}$, which matches the $L^2(\pi_{\theta_0})$-integrability threshold for $x^{-1}$. Their asymptotic normality result also uses a polynomial sampling scheme with $T_n\Delta_n \to 0$. The identification of their asymptotic covariance, and of the drift block in Cheng et al. (2024), with $I_0^{-1}$ is discussed in Section 1 of the main article. Prykhodko and Ralchenko (2025) establish consistency of discrete drift estimators for all positive drift parameters, and joint asymptotic normality under a polynomial sampling scheme with $T_n\Delta_n \to 0$.

Mazzonetto and Nieto (2025) study drift estimation for threshold CKLS processes from continuous and discrete observations. Under stationary initialization, suitable moment assumptions, and $T_n\Delta_n \to 0$, they obtain asymptotic normality of discrete estimators and LAN-type expansions of discretized likelihood and quasi-likelihood criteria. These LAN-type expansions concern discretized continuous-observation criteria, whereas Theorem 1.1 concerns the exact sampled likelihood.